\documentclass[11pt]{article}
\usepackage{amsbsy}
\usepackage{tikz}
\usepackage{amsmath, amssymb, amsfonts, euscript, pifont, mathrsfs, latexsym, graphicx, txfonts, cite, yhmath, url, color, relsize, float}
\usepackage[labelsep=period, font=footnotesize, labelfont=bf]{caption}

\usepackage[a4paper, left=2.5cm, right=2.5cm, top= 3cm, bottom=3cm]{geometry}

\newtheorem{construction}{Construction}[section]
\newtheorem{theorem}[construction]{Theorem}

\newtheorem{conjecture} {Conjecture}
\newtheorem{corollary} [construction]{Corollary}

\newtheorem{lemma}[construction] {Lemma}

\DeclareMathAlphabet{\mathpzc}{OT1}{pzc}{m}{it}

\newcommand{\mnull}[1]{\operatorname{\mathpzc{null}}(#1)}

\newcommand{\tr}[1]{\operatorname{\mathpzc{tr}}(#1)}

\begin{document}
\title{\vspace*{3cm} Short note on Randi\'{c} energy}
\author{Luiz Emilio Allem{\footnotesize$^\dagger$} \and Gonzalo Molina{\footnotesize$^\ddagger$} \and Adri\'{a}n Pastine{\footnotesize$^\ddagger$}}
\date{}
\maketitle
\vspace*{-8mm}
\begin{center}

{\footnotesize{ $\dagger$ Instituto de Matem\'atica, UFRGS, Porto Alegre, RS, 91509--900, Brazil}
\vspace{2mm}

{ $^\ddagger$Universidad Nacional de San Luis, Departamento de Matem\'aticas, San Luis, Argentina}
\vspace{2mm}

{\tt adrian.pastine.tag@gmail.com},~~~ {\tt emilio.allem@ufrgs.br}, ~~~{\tt lgmolina@unsl.edu.ar}}
\end{center}
\vspace{6mm}


\baselineskip = .30in

\begin{abstract}

In this paper, we consider the Randi\'{c} energy $RE$ of simple connected graphs.
We provide upper bounds for $RE$ in terms of the number of vertices and the nullity of the graph.
We present families of graphs that satisfy the Conjecture proposed by Gutman, Furtula and Bozkurt \cite{Gutman1} about the maximal RE. For example, we show that starlikes of odd order satisfy the conjecture.
\end{abstract}



\section{Introduction}
The problem of finding the graphs with maximal and minimal energy has been extensively studied for several matrices. For the Adjacency matrix, Gutman \cite{Gutman2} proved the following.
\begin{theorem}
Let $T$ be a tree on $n$ vertices. Then
$$E(S_{n})\leq E(T)\leq E(P_{n}).$$
\end{theorem}
Where $P_{n}$ and $S_{n}$ stand for the $n$-vertex path and the $n$-vertex star.
Radenkovi\'{c} and Gutman \cite{Gutman3} conjecturated the following about the Laplacian energy.
\begin{conjecture}
Let $T$ be a tree on $n$ vertices. Then
$$LE(P_{n})\leq LE(T)\leq LE(S_{n}).$$
\end{conjecture}
Fritscher et al. \cite{Eliseu1} proved that among the trees the star has maximal Laplacian energy. The problem of minimal Laplacian energy is still open. 

In the paper of Gutman, Furtula and Bozkurt \cite{Gutman1} on the energy of the Randi\'{c} matrix, graphs called sun and double sun were defined. The authors presented the following conjecture about the connected graphs with maximal RE.
\begin{conjecture}
\label{conjecturerandic}
Let $G$ be a connected graph on $n$ vertices. Then
  \begin{equation*}
RE(G) \leq \begin{cases}
RE(Sun)&\text{if $n$ is odd,}\\
RE(Balanced\mbox{ }double\mbox{ }sun)&\text{if $n$ is even.}
\end{cases}
\end{equation*}
\end{conjecture}

In this article we present bounds for the Randi\'{c} energy. And some families of graphs that satisfy the conjecture above.

Let $G=(V,E)$ be an undirected graph with vertex set $V$ and edge set $E$.
The Randi\'{c} matrix $R = [r_{ij}]$ of a graph $G$ is defined
\cite{Bozkurt2010a,Das2016,Gutman1}
as
$$
r_{ij} =
\displaystyle{
        \left\{
               \begin{array}{cl}
                  \frac{1}{\sqrt{d_{u}d_{v}}}  &  \mbox{ if } uv\in E \\
                            0                  &  \mbox{ otherwise }
               \end{array}
       \right.}
$$
Denote the eigenvalues of $R$ by $\lambda_{1}, \ldots, \lambda_{n}$. The multiset $\sigma_{R}=\{\lambda_{1}, \ldots, \lambda_{n}\}$ will be called the $R$-spectrum of the graph $G$.

The  Randi\'{c} energy $RE(G)$ of a graph $G$ is
$$\sum_{i=1}^{n} |\lambda_{i}|.$$

Historically, the $RE$ is related to a
descriptor for molecular graphs used by Milan Randi\'{c} in 1975 \cite{Randic1975}.
The normalized Laplacian matrix, defined by Chung \cite{Chung1997}, can be written using the Randi\'{c} matrix as
$$\mathcal{L}=I _{n}-R.$$
And the eigenvalues of $\mathcal{L}$ are given by
$$\mu_{i}=1-\lambda_{i}$$
for $i=1\ldots n$. For graphs without isolated vertices Cavers \cite{Cavers} defined the normalized Laplacian energy as
$$E_{\mathcal{L}}(G)=\sum_{i=1}^{n}|\mu_{i}-1|.$$
An interesting fact about $E_{\mathcal{L}}(G)$, see \cite{Gutman1}, is that if $G$ does not have isolated vertices then $$RE(G) = E_{\mathcal{L}}(G).$$
Thus, results in this paper on Randi\'{c} energy apply also to normalized Laplacian energy.

This paper is organized as follows. In Section \ref{Sun}, we present closed formulas for the Randi\'{c} energy of the sun and the double sun. In Section \ref{Spectrum}, we use some known eigenvalues to  provide upper bounds for $RE$ in terms of the number of vertices and the nullity of the graph. In Section \ref{bounds}, we use bounds for the Randi\'{c} index $R_{-1}(G)$ to improve bounds for $RE$. In Section \ref{TB}, we show that some families of graphs, for example starlikes of odd order, satisfy the conjecture proposed in \cite{Gutman1}.

\section{Randi\'{c} energy of sun and double sun}
\label{Sun}
In the work of Gutman, Furtula and Bozkurt \cite{Gutman1}
on the energy of the Randi\'{c} matrix, two families of trees were defined, sun and doble sun.
For each $p \geq 0$, the $p$-{\em sun},
which we denote with $S^p$,
is the tree of order $n = 2p + 1$
formed by taking the star on $p+1$ vertices and subdividing each edge.
\begin{figure}[h]
\centering
\begin{tikzpicture}
\path( 0,0)node[shape=circle,draw,fill=black] (primeiro) {}
      ( 1,0.5)node[shape=circle,draw,fill=black] (segundo) {}
      ( 1,1)node[shape=circle,draw,fill=black] (terceiro) {}
      (2,0.5)node[shape=circle,draw,fill=black] (quarto){}
      (2,1)node[shape=circle,draw,fill=black] (quinto) {}
      ( 1,-0.5)node[shape=circle,draw,fill=black] (sexto){}
      ( 2,-0.5)node[shape=circle,draw,fill=black] (setimo){};

      \draw[-](primeiro)--(segundo);
      \draw[-](primeiro)--(terceiro);
      \draw[-](primeiro)--(sexto);
      \draw[-](terceiro)--(quinto);
      \draw[-](segundo)--(quarto);
      \draw[-](setimo)--(sexto);

			\draw[loosely dotted,line width=0.5mm] (1.5,0.4) -- (1.5,-0.4);
\end{tikzpicture}
\caption{\label{fig-sun} Sun.}
						
\end{figure}

For $p, q \geq 0$ the $(p,q)$-{\em double sun},
denoted $D^{p,q}$,
is the tree of order $n = 2(p + q + 1)$ obtained by connecting the
centers of $S^p$ and $S^q$ with an edge.
Without loss of generality we assume $p \geq q$.
When $p-q \leq 1$ the double sun is called {\em balanced}.
\begin{figure}[h]
\centering
\begin{tikzpicture}
\path( 0,0)node[shape=circle,draw,fill=black] (primeiro) {}
      ( 1,0.5)node[shape=circle,draw,fill=black] (segundo) {}
      ( 1,1)node[shape=circle,draw,fill=black] (terceiro) {}
      (2,0.5)node[shape=circle,draw,fill=black] (quarto){}
      (2,1)node[shape=circle,draw,fill=black] (quinto) {}
      ( 1,-0.5)node[shape=circle,draw,fill=black] (sexto){}
      ( 2,-0.5)node[shape=circle,draw,fill=black] (setimo){}
     ( -1,0)node[shape=circle,draw,fill=black] (oitavo) {}
      ( -2,0.5)node[shape=circle,draw,fill=black] (nono) {}
      ( -2,1)node[shape=circle,draw,fill=black] (decimo) {}
      (-3,1)node[shape=circle,draw,fill=black] (decprimeiro){}
      (-3,0.5)node[shape=circle,draw,fill=black] (decseg) {}
      ( -3,-0.5)node[shape=circle,draw,fill=black] (decterc){}
      ( -2,-0.5)node[shape=circle,draw,fill=black] (decquar){};

       \draw[-](primeiro)--(oitavo);
      \draw[-](primeiro)--(segundo);
      \draw[-](primeiro)--(terceiro);
      \draw[-](primeiro)--(sexto);
      \draw[-](terceiro)--(quinto);
      \draw[-](segundo)--(quarto);
      \draw[-](setimo)--(sexto);
			\draw[-](oitavo)--(nono);
      \draw[-](oitavo)--(decimo);
      \draw[-](oitavo)--(decquar);
      \draw[-](decimo)--(decprimeiro);
      \draw[-](nono)--(decseg);
      \draw[-](decterc)--(decquar);

			\draw[loosely dotted,line width=0.5mm] (1.5,0.4) -- (1.5,-0.4);
			\draw[loosely dotted,line width=0.5mm] (-2.5,0.4) -- (-2.5,-0.4);
			\end{tikzpicture}
						\caption{\label{fig-doublesun} Double Sun.}
						
\end{figure}
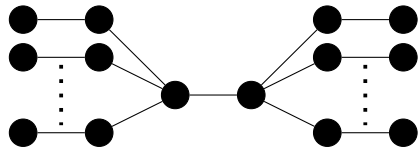

In \cite{Gutman1} was conjectured that the connected graph with maximal Randi\'{c} energy is a tree. More specifically, if $n\geq 1$ is odd, the sun is conjectured to have greatest Randi\'{c} energy among graphs with $n$ vertices. And, if $n\geq 2$ is even, then the balanced double sun is conjectured to have greatest Randi\'{c} energy among graphs with $n$ vertices.

Using the algorithm developed in \cite{Braga1}, for locating eigenvalues in trees for the normalized Laplacian matrix, we can compute the characteristic polynomials of the sun and the balanced double sun.

The characteristic polynomial of the sun with $p\geq 1$ is:

$$\mbox{det}(\lambda I-\mathcal{L})=(-1)(\lambda-(\frac{2+\sqrt{2}}{2}))^{p-1}(\lambda-(\frac{2-\sqrt{2}}{2}))^{p-1}(\lambda)(\lambda-2)(\lambda-1).$$

It follows that

$$E_{\mathcal{L}}(S^{p})=\sum_{i=1}^{n}|\lambda_{i}(l)-1|=2(p-1)\frac{\sqrt{2}}{2}+2=(n-3)\frac{\sqrt{2}}{2}+2.$$

Suppose that $p\geq q$ and $p+q\geq 2$. Then, the characteristic polynomial of  $D^{p,q}$ is:
$$\mbox{det}(\lambda I-\mathcal{L})=\lambda(\lambda-2)(\lambda-(\frac{2+\sqrt{2}}{2}))^{p+q-2}(x-(\frac{2-\sqrt{2}}{2}))^{p+q-2}(q(\lambda))$$
with $$q(\lambda)=\lambda^{4}-4\lambda^{3}+\frac{1}{4}\frac{(22p+20qp+22q+20)}{(q+1)(p+1)}\lambda^{2}+\frac{1}{4}\frac{(-12p-8qp-12q-8)}{(q+1)(p+1)}+\frac{1}{4}\frac{(1+2p+2q)}{(q+1)(p+1)}.$$
 It is known that the graph $G$ is bipartite if and only if for each normalized laplacian eigenvalue $\lambda$, the value $2-\lambda$ is also an eigenvalue of $G$. Using this fact, we can write $q(\lambda)$ as 
 $$q(\lambda)=(\lambda-\lambda_{a})(\lambda-\lambda_{b})(\lambda-(2-\lambda_{a}))(\lambda-(2-\lambda_{b}))$$ with $\lambda_{a}\leq\lambda_{b}$. Now, we can compute the energy of the balanced double sun in both cases as follows:

If $p=q$ then $$E_{\mathcal{L}}(D^{p,p})=\frac{\sqrt{2}(n^{2}-4n-12)+4\sqrt{n^{2}+4n+20}}{2(n+2)}.$$

If $q=p-1$ then
\begin{align*}
E_{\mathcal{L}}(D^{p,p-1})=\frac{\sqrt{2}}{2n(n+4)}(n^{3}-2n^{2}-24+2\sqrt{n(n+4)(n^{2}+8+\sqrt{-64n+n^{4}+64})}+ \\
2\sqrt{n(n+4)(n^{2}+8-\sqrt{-64n+n^{4}+64})}).
\end{align*}

Now, we can rewrite the Conjecture \ref{conjecturerandic} for the Randi\'{c} energy using closed formulas.

\begin{conjecture}
Let $G$ be a connected graph on $n$ vertices. Then for $k\geq 3$ odd we have that
\begin{equation*}
RE(G) \leq \begin{cases}
RE(S^{p})= E_{\mathcal{L}}(S^{p})&\text{if $n=k$}\\
RE(DS^{p,p})= E_{\mathcal{L}}(D^{p,p-1})&\text{if $n=2k$}\\
RE(DS^{p,p-1})=E_{\mathcal{L}}(D^{p,p-1}) &\text{if $n=2k+2$}
\end{cases}
\end{equation*}
\end{conjecture}
\section{Upper bounds for $RE$}
\label{Spectrum}
In this section, we present upper bounds for $RE$ in terms of the number of vertices and the nullity of $G$. The main tool we use to study the Randi\'{c} energy of graphs is the trace of $R^2$, taking advantage of the eigenvalues
we know for $G$.

The next Theorem is  the generalized mean inequality that will be used in the results following it.
\begin{theorem}
\label{generalizedmean}
If $x_1,\ldots,x_n$ are nonnegative real numbers, $p$ and $q$ are positive integers, and $p<q$, then
\[
\left(\frac{1}{n}\sum_{i=1}^n x_i^p\right)^{1/p}\leq \left(\frac{1}{n}\sum_{i=1}^n x_i^q\right)^{1/q}
\]
\end{theorem}

The next Lemma can also be found in \cite{Cavers}. But for completeness we present a proof here.

\begin{lemma}
\label{oldtrbound}
Let $G$ be a graph of order $n$ with no isolated vertices. Then
$$RE(G)\leq \sqrt{n\tr{R^2}}.$$
\end{lemma}
\begin{Prf}
Applying Theorem \ref{generalizedmean} with $p=1$, $q=2$, $x_i=|\lambda_i|$ yields
\begin{align*}
\frac{1}{n}\sum_{i=1}^n |\lambda_i|&\leq \left(\frac{1}{n}\sum_{i=1}^n |\lambda_i|^2\right)^{1/2}\\
\sum_{i=1}^n |\lambda_i|&\leq n\frac{1}{\sqrt{n}}\left(\sum_{i=1}^n |\lambda_i|^2\right)^{1/2}\\
\sum_{i=1}^n |\lambda_i|&\leq \sqrt{n}\left(\sum_{i=1}^n |\lambda_i|^2\right)^{1/2}.
\end{align*}
And the last inequality is exactly $RE(G)\leq \sqrt{n\tr{R^2}}$.
\end{Prf} $\Box$

Notice that Lemma \ref{oldtrbound} can be improved when some of the eigenvalues for $R$ are known.
Consider $\Psi$ a sub-multiset of $\sigma_R$ and denote the multiset difference by $\sigma_R\setminus\Psi$.
\begin{theorem}\label{generaltrbound}
Let $G$ be a graph, and let $\Psi$ be a sub-multiset of $\sigma_R$. Then
\[
RE(G)\leq \sqrt{(n-|\Psi|)\left(\tr{R^2}-\sum_{\lambda\in\Psi}\lambda^2\right)}+\sum_{\lambda\in\Psi}|\lambda|.
\]
\end{theorem}
\begin{Prf}
Notice that
\[
RE(G)=\sum_{i=1}^n |\lambda_i|=\sum_{\lambda\in\sigma_R\setminus\Psi}|\lambda|+\sum_{\lambda\in\Psi}|\lambda|.
\]
Applying Theorem \ref{generalizedmean} to the elements of $\sigma_R\setminus\Psi$ yields
\begin{align*}
\frac{1}{|\sigma_R\setminus\Psi|}\sum_{\lambda\in\sigma_R\setminus\Psi}|\lambda|\leq & \left( \frac{1}{|\sigma_R\setminus\Psi|}\sum_{\lambda\in\sigma_R\setminus\Psi}(|\lambda|)^2\right)^{1/2}\\
\sum_{\lambda\in\sigma_R\setminus\Psi}|\lambda|\leq & \left( |\sigma_R\setminus\Psi|\sum_{\lambda\in\sigma_R\setminus\Psi}(\lambda)^2\right)^{1/2}.\\
\end{align*}
But $|\sigma_R\setminus\Psi|=n-|\Psi|$, and
$\sum_{\lambda\in\sigma_R\setminus\Psi}(\lambda)^2=\tr{R^2}-\sum_{\lambda\in\Psi}\lambda^2$. Thus
\[
\sum_{\lambda\in\sigma_R\setminus\Psi}|\lambda|\leq \left( (n-|\Psi|)\left(\tr{R^2}-\sum_{\lambda\in\Psi}\lambda^2\right)\right)^{1/2}.
\]
Hence
\begin{align*}
RE(G)\leq \sqrt{(n-|\Psi|)\left(\tr{R^2}-\sum_{\lambda\in\Psi}\lambda^2\right)}+\sum_{\lambda\in\Psi}|\lambda|.
\end{align*}
\end{Prf} $\Box$

We can apply Theorem \ref{generaltrbound} to graphs in general, using that $1$ is an eigenvalues of $R$ for every graph $G$, and using that $-1$ is an eigenvalue whenever the graph is bipartite. Furthermore, we can use the dimension of the null space, denoted by $\mnull{R}$, as that counts the multiplicity
of $0$ as an eigenvalue. Notice that $\sum_{\lambda\in\Psi}\lambda^2=\sum_{\lambda\in\Psi}|\lambda|=1$
in the general case, and  $\sum_{\lambda\in\Psi}\lambda^2=\sum_{\lambda\in\Psi}|\lambda|=2$ in the bipartite
case.
Hence, we obtain the following.
\begin{corollary}\label{trspecificbound}
If $G$ is a graph, then
\[
RE(G)\leq \sqrt{(n-1-\mnull{R})(\tr{R^2}-1)}+1.
\]
Furthermore, if $G$ is bipartite, then
\[
RE(G)\leq \sqrt{(n-2-\mnull{R})(\tr{R^2}-2)}+2.
\]
\end{corollary}
Corollary \ref{trspecificbound} will be our main tool to bound $RE(G)$ in the following sections.

\section{Upper bounds using $R_{-1}(G)$}
\label{bounds}
The randi\'{c} index of $G$, denoted by $R_{-1}(G)$, satisfies the equality
\[
R_{-1}(G)=\frac{1}{2}\tr{R^2}.
\]
Hence, any upper bound of $R_{-1}(G)$ may yield an upper bound for $RE(G)$.

In the next Theorem we summarize some upper bounds for $R_{-1}(G)$.

\begin{theorem}
\label{boundonrandicindex}
In \cite{Cavers}, Cavers et al. showed that if $G$ is a connected graph on $n\geq 3$ vertices. Then
\[
R_{-1}(G)\leq \frac{15(n+1)}{56}.
\]

In \cite{Clark}, Clark and Moon proved that if $T$ is a tree, then
\[R_{-1}(T)\leq \frac{5n+8}{18}.\]

In \cite{Li,Pav}, they proved that if $T$ is a tree of order $n\geq 103$, then
\[
R_{-1}(T)\leq \frac{15n-1}{56}.\]
\end{theorem}

Applying Corollary \ref{trspecificbound} together with Theorem \ref{boundonrandicindex} yields
\begin{corollary}\label{boundsforallgraphs}
Let $G$ be a graph, then
\[
RE(G)\leq \sqrt{(n-1-\mnull{R})\frac{15n-13}{28}}+1.
\]

Let $G$ be a bipartite graph, then
\[
RE(G)\leq \sqrt{(n-2-\mnull{R})\frac{15n-41}{28}}+2.
\]

Let $T$ be a tree, then
\[
RE(T)\leq \sqrt{(n-2-\mnull{R})\frac{5n-10}{9}}+2.
\]

Let $T$ be a tree of order $n\geq 103$, then
\[
RE(T)\leq \sqrt{(n-2-\mnull{R})\frac{15n-57}{28}}+2.
\]
\end{corollary}

It is known that trees of odd order have nullity at least one.
If we consider $\mnull{R}=0$ for trees of even order and $\mnull{R}=1$ for trees of odd order in  Corollary \ref{boundsforallgraphs} we obtain the following result.
\begin{corollary}
\label{null}
Let $T$ be a tree of even order $n\geq 2$, then
\begin{equation}
\label{eq1}
RE(T)\leq \sqrt{(n-2)\frac{5n-10}{9}}+2.
\end{equation}

Let $T$ be a tree of odd order $n\geq 2$, then
\begin{equation}
\label{eq1nul1}
RE(T)\leq \sqrt{(n-3)\frac{5n-10}{9}}+2.
\end{equation}

Let $T$ be a tree of even order $n\geq 103$, then 
\begin{equation}
\label{eq2}
RE(T)\leq \sqrt{(n-2)\frac{15n-57}{28}}+2.
\end{equation}

Let $T$ be a tree of odd order $n\geq 103$, then 
\begin{equation}
\label{eq2nul1}
RE(T)\leq \sqrt{(n-3)\frac{15n-57}{28}}+2.
\end{equation}
\end{corollary}

In \cite{Das2016} the following bound for the energy of trees was found.
\begin{theorem}\cite{Das2016}
\label{TheoremDas}
Let $T$ be a tree of order $n$, then
\begin{equation}
\label{eq3}
RE(T)\leq 2\sqrt{\left\lfloor\frac{n}{2}\right\rfloor\frac{5n+8}{18}}.
\end{equation}
\end{theorem}

With a simple calculation we can compare the bounds in Corollary \ref{null} and Theorem \ref{TheoremDas}.
If $n$ is even we use $\left\lfloor n/2\right\rfloor=n/2$. Then $(\ref{eq1})<(\ref{eq3})$ if $n\geq 3$ and $(\ref{eq1})=(\ref{eq3})$ if $n=2$.
If $n$ is odd we use $\left\lfloor n/2\right\rfloor=(n-1)/2$. Then $(\ref{eq1nul1})<(\ref{eq3})$ if $n\geq 3$.

Corollary \ref{boundsforallgraphs} also suggests that if $\mnull{R}$ is sufficiently large, then
$RE(G)\leq (n-3)\frac{\sqrt{2}}{2}+2$. The following theorem gives a lower bound for the
nullity that yields the inequality.
\begin{theorem}\label{nullbound}
If $\mnull{R}\geq \dfrac{n^2+56n-141-28(n-3)\sqrt{2}}{15n-13}$, then $RE(G)\leq (n-3)\frac{\sqrt{2}}{2}+2$.
\end{theorem}
\begin{Prf}
If $G$ is a graph, then writing Corollary \ref{trspecificbound} in terms of $R_{-1}(G)$ yields
\[
RE(G)\leq \sqrt{n-1-\mnull{R}}\sqrt{2R_{-1}(G)-1}+1.
\]

Using Theorem \ref{boundonrandicindex},
\begin{align*}
RE(G)\leq& \sqrt{n-1-\mnull{R}}\sqrt{\frac{15(n+1)}{28}-1}+1\\
=&\sqrt{n-1-\mnull{R}}\sqrt{\frac{15n-13}{28}}+1.
\end{align*}
Thus, we want
\begin{align*}
\sqrt{n-1-\mnull{R}}\sqrt{\frac{n-13}{28}}+1&\leq (n-3)\frac{\sqrt{2}}{2}+2,\\
(n-1-\mnull{R})\frac{15n-13}{28}&\leq ((n-3)\frac{\sqrt{2}}{2}+1)^2,\\
(n-1-\mnull{R})\frac{15n-13}{28}&\leq (n^2-6n+9)\frac{1}{2}+(n-3)\sqrt{2}+1,\\
(n-1-\mnull{R})(15n-13)&\leq 14 n^2 - 84n + 126+28(n-3)\sqrt{2}+28,\\
(n-1-\mnull{R})&\leq \frac{14 n^2 - 84n + 154+28(n-3)\sqrt{2}}{15n-13}.\\
\end{align*}
Hence,
\begin{align*}
-\mnull{R}&\leq \frac{14 n^2 - 84n + 154+28(n-3)\sqrt{2}}{15n-13}-n+1,\\
\mnull{R}&\geq -\frac{14 n^2 - 84n + 154+28(n-3)\sqrt{2}}{15n-13}+n-1,\\
\mnull{R}&\geq \frac{15n^2-13n-15n+13-14n^2+84n-154-28(n-3)\sqrt{2}}{15n-13},\\
\mnull{R}&\geq \frac{n^2+56n-141-28(n-3)\sqrt{2}}{15n-13}.\\
\end{align*}
\end{Prf} $\Box$

In the case of trees, using $R_{-1}(G)\leq (15n-1)/56$ and the fact that $\pm 1$ are eigenvalues,
Theorem \ref{nullbound} can be improved to show that $\mnull{T} \geq \dfrac{n^2-3n-12}{15n-57}$ implies $RE(T)\leq (n-3)\dfrac{\sqrt{2}}{2}+2$.


A \textit{suspended path} is a path $uvw$, with $d_u=1$ and $d_v=2$, i.e. $u$ is a pendent vertex and its neighbor
has degree $2$.
The next result improves the bound on $R_{-1}(G)$ when $G$ has no suspended paths.

\begin{theorem}\cite{Cavers}
Let $G$ be a connected graph on $n\geq 3$ vertices. If $G$ has no suspended paths,
then
\[
R_{-1}(G)\leq \dfrac{n}{4}
\]
\end{theorem}

If $G$ is bipartite, then we can use that $\pm 1$ are eigenvalues of $R$, and, hence, $1$ is an eigenvalue of multiplicity $2$
of $R^2$ to obtain the following.
\begin{theorem}
Let $G$ be a connected bipartite graph on $n\geq 3$ vertices. If $G$ has no suspended paths, then
\[
RE(G)\leq \sqrt{n-2}\sqrt{n-4}\dfrac{\sqrt{2}}{2}+2.
\]
\end{theorem}
Notice that $\sqrt{n-2}\sqrt{n-4}=\sqrt{n^2-6n+8}< \sqrt{n^2-6n+9}=(n-3)$.
Therefore, if $n=2p+1$ is odd and $G$ is a connected bipartite graph on $n\geq 3$ vertices that has
no suspended paths, then
$RE(G)\leq RE(S^{p})$.

When $G$ is not bipartite, it is better to look at a result using $\mnull{R}$.
\begin{theorem}
Let $G$ be a connected graph on $n\geq 3$ vertices. If $G$ has no suspended paths, then
\[
RE(G)\leq \sqrt{n-1-\mnull{R}}\sqrt{n-2}\dfrac{\sqrt{2}}{2}+1.
\]
\end{theorem}
Notice in particular that if $\mnull{R}\geq 1$, then
\begin{align*}
\sqrt{n-1-\mnull{R}}\sqrt{n-2}\dfrac{\sqrt{2}}{2}+1\leq &\sqrt{n-2}\sqrt{n-2}\frac{\sqrt{2}}{2}+1\\
= & (n-2)\frac{\sqrt{2}}{2}+1\\
= & (n-3)\frac{\sqrt{2}}{2}+\frac{\sqrt{2}}{2}+1\\
< (n-3)\frac{\sqrt{2}}{2}+2.\\
\end{align*}
Hence, if $n=2p+1$ is odd and $G$ is a connected graph on $n\geq 3$ vertices that has no suspended paths,
and $\mnull{R}\geq 1$, then $RE(G)<RE(S^{p})$.
\section{TB graphs}
\label{TB}
In the previous section we showed how finding good bounds for $\tr{R^2}$ yields good bounds for $RE(G)$.
In this section we study $\tr{R^2}$ for a particular family of bipartite graphs, and use it to show
that their randi\'{c} energy is bounded by the randi\'{c} energy of the sun graph.
The family that we consider is bipartite graphs with bipartition $A,B$, such that $\deg(b)\leq 2$ for every $b\in B$.
We denote this graphs as TB graphs.
Notice that the family of TB graphs include many important subfamilies of graphs:

\begin{itemize}
\item starlike trees, which are trees with exactly one vertex of degree greater than 2;
\item basic trees, see \cite{JMS}, which are trees with a unique maximum independent set of size $\lceil n/2 \rceil$
(the importance of this trees is due to their null space);
\item a graph obtained by taking a graph $G$ and replacing each edge $e=\{v_1,v_2\}$ by two edges $e_1=\{v_1,w_e\}$ and
$e_2=\{v_2,w_e\}$.
\end{itemize}

The graphs described above satisfy the condition $\deg(b)=2$ for every $b\in B$.
In this section, we give a bound on $\tr{R^2}$ for any TB graph. Before doing so, we
give a short explanation on how to find $\tr{R^2}$ for any bipartite graph.
Let $G$ be a bipartite graph.
As $G$ is bipartite, the underlying graph of $R^2$ has two connected components.

Consider a bipartite graph $G$ with $V(G)=A\cup B$. Let $R$ be the Randi\'{c} matrix of $G$ indexed first with the vertices in $A$ and then in $B$. As there are no edges between vertices in $A$ and no edges between
vertices in $B$, $R$ is a block anti-diagonal matrix. I.e., $R$ is of the form
$$R=\begin{bmatrix}
0 & C \\
C^{t} & 0
\end{bmatrix},$$
where $C$ is a $|A|\times |B|$ matrix.
Then
$$R^{2}=\begin{bmatrix}
R_{A}^{2} & 0 \\
0 & R_{B}^{2}
\end{bmatrix}$$
with $R_{A}^{2}=CC^{t}$ and $R_{B}^{2}=C^{t}C$. It follows that $tr(R_{A}^{2})=tr(R_{B}^{2})$, and  
$tr(R^{2})=2tr(R_{A}^{2})=2tr(R_{B}^{2})$.
There is actually
a big difference between vertices of degree $1$ and vertices of degree $2$ in $B$, hence we partition $B$ into
$B_1=\{b\in B\,|\,\deg(b)=1\}$ and $B_2=\{b\in B\,|\, \deg(b)=2\}$.

\begin{lemma}\label{lemmaN(a)}
Let $G$ be a connected TB graph with $|G|\geq 3$. Then for every $a\in A$,
\[
\frac{1}{2}\leq R^2_{a,a}\leq\frac{1}{2}+\frac{1}{4}|N(a)\cap B_1|,
\]
where $N(a)$ is the neighborhood of $a$.
\end{lemma}
\begin{Prf}
If $\deg(a)\geq 2$, then
\begin{align*}
 R^2_{a,a}=&\sum_{b\in N(a)}\frac{1}{\deg(b)\deg(a)}\\
 =&\sum_{b\in N(a)\cap B_1}\frac{1}{\deg(a)}+\sum_{b\in N(a)\cap B_2}\frac{1}{2\deg(a)}\\
 =& \sum_{b\in N(a)\cap B_1}\frac{2}{2\deg(a)}+\sum_{b\in N(a)\cap B_2}\frac{1}{2\deg(a)}\\
  =& \sum_{b\in N(a)\cap B_1}\frac{1}{2\deg(a)}+\sum_{b\in N(a)\cap B_1}\frac{1}{2\deg(a)}+\sum_{b\in N(a)\cap B_2}\frac{1}{2\deg(a)}\\
=& \sum_{b\in N(a)\cap B_1}\frac{1}{2\deg(a)}+\sum_{b\in N(a)}\frac{1}{2\deg(a)}\\
=& |N(a)\cap B_1|\frac{1}{2\deg(a)}+\deg(a)\frac{1}{2\deg(a)}\\
=& |N(a)\cap B_1|\frac{1}{2\deg(a)}+\frac{1}{2},\\
\end{align*}
thus
\[
\frac{1}{2}\leq R^2_{a,a}\leq |N(a)\cap B_1|\frac{1}{4}+\frac{1}{2},
\]
where the second inequality follows from  $\deg(a)\geq 2$.

If $\deg(a)=1$, let $b$ be the only neighbor of $a$,
\begin{align*}
R^2_{a,a}=&\frac{1}{\deg(b)}\\
=& \frac{1}{2}\\
=& |N(a)\cap B_1|\frac{1}{4}+\frac{1}{2},
\end{align*}
because $\deg(b)= 2$, as $G$ is a connected TB graph with at least 3 vertices.
\end{Prf} $\Box$

Notice that if $G$ is a TB graph, then $E(G)=|B_1|+2|B_2|$. As $G$ is connected,
$E(G)\geq n-1=|A|+|B_1|+|B_2|-1$. Thus $2|B_2|\geq |A|+|B_2|-1$, or $|A|\leq |B_2|+1$.
We can now bound the trace.
\begin{lemma}\label{lemmatraceR2}
Let $G$ be a connected TB graph with $|G|\geq 3$. Then $\tr{R_A^2}\leq \frac{n+1}{4}$.
\end{lemma}
\begin{Prf}
\[
\tr{R_A^2}=\sum_{a\in A}R_{a,a}^2
\leq \sum_{a\in A}\left(\frac{1}{2}+\frac{1}{4}|N(a)\cap B_1|\right)
\]
but as the vertices in $B_1$ have degree $1$, they each appear in exactly one $N(a)\cap B_1$. Hence
\[
\sum_{a\in A}\frac{1}{4}|N(a)\cap B_1|=\sum_{b\in B_1}\frac{1}{4}=\frac{|B_1|}{4}.
\]

Then
\begin{align*}
\tr{R_A^2}\leq&  \sum_{a\in A}\left(\frac{1}{2}\right)+\frac{1}{4}|B_1|\\
\leq& \frac{2|A|}{4}+\frac{1}{4}|B_1|\\
\leq& \frac{|A|+|B_2|+1}{4}+\frac{1}{4}|B_1|\\
\leq& \frac{|A|+|B_2|+1+| B_1|}{4}\\
\leq& \frac{n+1}{4}.
\end{align*}
\end{Prf} $\Box$

As a TB graph is bipartite, $\tr{R^2}=2\tr{R^2_A}$.
\begin{lemma}\label{cotatrTB}
Let $G$ be a connected TB graph, then $\tr{R^2}\leq \frac{n+1}{2}$.
\end{lemma}

\begin{lemma}\label{lemmacotaconnnul}
Let $G$ be a bipartite graph. If $\tr{R^2}\leq \frac{n+1}{2}$, then
\[
RE(G)\leq \sqrt{n-2-\mnull{R}}\sqrt{n-3}\frac{\sqrt{2}}{2}+2.\]
\end{lemma}
\begin{Prf}
As $G$ is bipartite, Corollary \ref{trspecificbound} yields
\begin{align*}
RE(G)\leq \sqrt{(n-2-\mnull{R})\tr{R^2-2}}+2.
\end{align*}

But $\tr{R^2}-2\leq \frac{n+1}{2}-2=\frac{n-3}{2}$.
Thus
\begin{align*}
RE(G)\leq \sqrt{n-2-\mnull{R}}\sqrt{n-3}\frac{\sqrt{2}}{2}+2.
\end{align*}
\end{Prf} $\Box$

We can now combine Lemma \ref{cotatrTB} with Lemma \ref{lemmacotaconnnul} to obtain the following.
\begin{theorem}
Let $G$ be a connected TB graph. Then
\[
RE(G)\leq \sqrt{n-2}\sqrt{n-3}\frac{\sqrt{2}}{2}+2.\]
Even more, if $\mnull{R}\geq 1$, then
\[RE(G)\leq (n-3)\frac{\sqrt{2}}{2}+2.\]
\end{theorem}

Notice that bipartite graphs with an odd number of vertices have $\mnull{R}\geq 1$. This shows that the randi\'{c}
energy of TB graphs of odd order is less or equal than the randi\'{c} energy of the sun graph of that same order.

%
%
%
%

\section*{Acknowledgments}
This work was partially supported by the Universidad Nacional de San Luis, grant PROICO 03-0918, and MATH AmSud, grant 18-MATH-01.

\end{document}